\documentclass[11pt]{article}
\usepackage{latexsym}
\usepackage{amssymb,amsmath,amsfonts}
\usepackage[mathscr]{eucal}
\usepackage{amscd}
\usepackage{amsfonts,cite}
\usepackage{epsfig, amsmath, amsthm, amssymb}
\usepackage{latexsym}
\usepackage{graphics,epstopdf}
\begin{document}
\thispagestyle{empty}
\parindent=0mm
\begin{center}
{\Large{\bf 2-iterated Sheffer polynomials$^{\star}$}}\\~~

{\bf Subuhi Khan{\footnote{Corresponding author;
E-mail:~subuhi2006@gmail.com (Subuhi Khan)\\ \hspace*{0.45cm} Second author; E-mail:~mumtazrst@gmail.com (Mumtaz Riyasat)}} and Mumtaz Riyasat}\\

Department of Mathematics\\
Aligarh Muslim University\\
Aligarh, India\\
\end{center}

\parindent=0mm

\noindent
{\bf Abstract:}~In this article, the 2-iterated Sheffer polynomials are introduced by means of generating function and operational representation.
Using the theory of Riordan arrays and relations between the Sheffer sequences and Riordan arrays, a determinantal definition for these polynomials
is established. The quasi-monomial and other properties of these polynomials are derived. The generating function, determinantal definition, quasi-monomial and other properties for some new members belonging to this family are also considered.\\

\vspace{.35cm}

\noindent
{\large{\bf 1.~Introduction and preliminaries}}
\vspace{.35cm}

\parindent=8mm

Sequences of polynomials play a fundamental role in mathematics. One of the most famous classes of polynomial sequences is the class of Sheffer sequences, which contains many important sequences such as those formed by Bernoulli polynomials, Euler polynomials, Abel polynomials, Hermite polynomials, Laguerre polynomials, {\em etc.} and contains the classes of associated sequences and Appell sequences as two subclasses. Roman {\em et. al.} in \cite{Roman, RomanRota} studied the Sheffer sequences systematically by the theory of modern umbral calculus.\\

Roman \cite{Roman1} further developed the theory of umbral calculus and generalized the concept of Sheffer sequences so that more special
polynomial sequences are included, such as the sequences related to Gegenbauer polynomials, Chebyshev polynomials and Jacobi polynomials.\\

Let $\mathbb{K}$ be a field of characteristic zero. Let $\mathcal{F}$ be the set of all formal power series in the variable $t$ over $\mathbb{K}$.
Thus an element of $\mathcal{F}$ has the form
$$f(t)=\sum\limits_{k=0}^\infty a_k t^k,\eqno(1.1)$$
where $a_k \in \mathbb{K}$ for all $k \in \mathbb{N}:=\{0,1,2,\ldots\}$. The order O$(f(t))$ of a power series $f(t)$ is the smallest integer $k$
for which the coefficient of $t^k$ does not vanish. The series $f(t)$ has a multiplicative inverse, denoted by ${f(t)}^{-1}$ or $\frac{1}{f(t)}$,
if and only if O$(f(t))=0$. Then $f(t)$ is called an invertible series. The series $f(t)$ has a compositional inverse, denoted by $\bar{f}(t)$
and satisfying $f(\bar{f}(t))=\bar{f}(f(t))=t$, if and only if  O$(f(t))=1$. Then $f(t)$ is called a delta series.\\

Let $f(t)$ be a delta series and $g(t)$ be an invertible series of the following forms:
$$f(t)=\sum\limits_{n=0}^\infty f_n \frac{t^n}{n!},~~~f_0=0,~f_1\neq0\eqno(1.2a)$$
and
$${g}(t)=\sum\limits_{n=0}^\infty g_n \frac{t^n}{n!},~~~g_0\neq0.\eqno(1.2b)$$

According to Roman \cite[p.18 (Theorem 2.3.4)]{Roman}, the polynomial sequence $s_n(x)$ is uniquely determined by two (formal) power series given by equations (1.2a) and (1.2b). The exponential generating function of $s_n(x)$ is then given by

$$\frac{1}{g(\bar{f}(t))}e^{x(\bar{f}(t))}=\sum\limits_{n=0}^\infty s_n(x)\frac{t^n}{n!},\eqno(1.3)$$
for all $x$ in $\mathbb{C}$, where $\bar{f}(t)$ is the compositional inverse of $f(t)$. The sequence $s_n(x)$ in equation (1.3) is the Sheffer sequence for the pair $({g}(t),f(t))$.\\

It has been given in \cite{HerSheff} that the Sheffer polynomials $s_n(x)$ are `quasi-monomial', for this see \cite{Steff,Dat}. The associated
raising and lowering operators are given by
$$\hat{M}_s=\left(x-\frac{{g'}(\partial_x)}{{g}(\partial_x)}\right)\frac{1}{f'(\partial_x)}\eqno(1.4a)$$
and
$$\hat{P}_s=f(\partial_x),\eqno(1.4b)$$
respectively, where $\partial_x:=\frac{\partial}{\partial x}$.\\

Also, for an invertible series $g(t)$ and a delta series $f(t)$, the sequence $(s_n(x))_{n \in \mathbb{N}}$ is Sheffer for the pair $(g(t),f(t))$
if and only if
$$\langle g(t){f(t)}^k |s_n(x) \rangle=c_n \delta_{n,k},\eqno(1.5)$$
for all $n,k \geq 0$, where $\delta_{n,k}$ is the kronecker delta. Particularly, the Sheffer sequence for $(1,f(t))$ is called the associated sequence
for $f(t)$ defined by the generating function of the form:
$$e^{x\bar{f}(t)}=\sum\limits_{n=0}^\infty \tilde{s}_n(x)\frac{t^n}{n!}\eqno(1.6)$$
and the Sheffer sequence for $(g(t),t)$ is called the Appell sequence \cite{Appell} for $g(t)$ defined by the generating function of the form \cite{Roman}:
$$\frac{1}{g(t)}e^{xt}=\sum\limits_{n=0}^\infty A_n(x)\frac{t^n}{n!}.\eqno(1.7)$$

Now, we recall some preliminaries from \cite{Wang}.\\

Let $P$ be the algebra of polynomials in the variable $x$ over $\mathbb{K}$ and $P^\star$ be the dual vector space of all functionals on $P$.
Let $(c_n)_{n \in \mathbb{N}}$ be a fixed sequence of nonzero constants. Then for each $f(t)$ in $\mathcal{F}$, we define a linear functional
$f(t)$ in $P^\star$ is defined by

$$\langle f(t)|x^n \rangle=c_n a_n\eqno(1.8)$$
and define a linear operator $f(t)$ on $P$ is defined by
$$f(t)x^n=\sum\limits_{k=0}^n \frac{c_n}{c_{n-k}}a_k x^{n-k}.\eqno(1.9)$$

Now, we give the concept of Riordan arrays, which was introduced by Shapiro {\em et. al.} \cite{Shapiro} and further studied by many authors, for example see \cite{Barry,Cheon,GouldTXHe,TXHe1,TXHe,Sprug1,Sprug2}.\\

For an invertible series $g(t)$ and a delta series $f(t)$, a generalized Riordan array with respect to the sequence $(c_n)_{n \in \mathbb{N}}$ is a pair $(g(t),f(t))$, which is an infinite, lower triangular array $(a_{n,k})_{0\leq k \leq n<\infty}$ according to the following rule:

$$a_{n,k}=\left[\frac{t^n}{c_n}\right]g(t)\frac{(f(t))^k}{c_k},\eqno(1.10)$$
where the functions $g(t)(f(t))^k/c_k$ are called the column generating functions of the Riordan array. Particularly, the classical Riordan arrays
corresponds to the case $c_n=1$ and the exponential Riordan arrays corresponds to the case of $c_n=n!$. One of the most important applications of the theory
of the Riordan arrays is to deal with summations of the form $\sum\limits_{k=0}^n a_{n,k}h_k$ \cite{Sprug1,Sprug2}.\\

Also, for any fixed sequence $(c_n)_{n \in \mathbb{N}}$, the set of all Riordan arrays $(g(t),f(t))$ is a group under
matrix multiplication and is called a Riordan group with respect to $(c_n)_{n \in \mathbb{N}}$. The identity of
this group is $(1,t)$ and the inverse of the array $(g(t),f(t))$ is $(1/g(\bar{f}(t)),\bar{f}(t))$.\\

Costabile {\em et. al.} in \cite{Cos1} proposed a determinantal definition for the classical Bernoulli polynomials.
Further, Costabile and Longo in \cite{Cos2} introduced the determinantal definition for the Appell sequences.
Later on, by using the theory of Riordan arrays and relation between the Sheffer sequences
and Riordan arrays, the determinantal definition of Appell sequences is extended to Sheffer sequences \cite{Wang}. The relations between the Sheffer sequences
and Riordan arrays for the case of classical Sheffer sequences and classical Riordan arrays are given in \cite{TXHe} and
for the case of generalized Sheffer sequences and generalized Riordan arrays are given in \cite{GouldTXHe,Wang1}. \\

Let $(s_n(x))_{n \in \mathbb{N}}$ be Sheffer for $(g(t),f(t))$ and suppose
$$x^n=\sum\limits_{k=0}^n a_{n,k}s_k(x),\eqno(1.11)$$
then $a_{n,k}$ is the $(n,k)$ entry of the Riordan array $(g(t),f(t))$ \cite{Wang1}.\\

In view of above fact, the following determinantal definition for the Sheffer sequences holds true \cite{Wang}:\\

Let $(s_n(x))_{n \in \mathbb{N}}$ be Sheffer for $(g(t),f(t))$, then we have
$$\begin{array}{l}
s_0(x)=\frac{1}{a_{0,0}},\hspace{5in}(1.12)\\
s_n(x)=\frac{{(-1)}^n}{a_{0,0}a_{1,1}\ldots a_{n,n}}\left|\begin{array}{cccccc}
 1  & x  & x^2  & \cdots & x^{n-1}   & x^n \\
\\
a_{0,0}  &  a_{1,0}  &  a_{2,0}  & \cdots & a_{n-1,0} &a_{n,0} \\
\\
 0  &  a_{1,1}  &  a_{2,1}  & \cdots & a_{n-1,1} &a_{n,1} \\
\\
 0 & 0 & a_{2,2} & \cdots & a_{n-1,2}&  a_{n,2}\\
 . & . & . & \cdots & . & . \\
 . & . & . & \cdots & . & . \\
 0 & 0 & 0 & \cdots & a_{n-1,n-1} & a_{n,n-1}
\end{array} \right|,\\
\\
\hspace{1.5cm}=\frac{{(-1)}^n}{a_{0,0}a_{1,1}\ldots a_{n,n}} \textrm{det}\left(\begin{array}{c}
                                                                       X_{n+1} \\
                                                                       S_{n \times (n+1)} \\
                                                                     \end{array}\right),\end{array}\eqno(1.13)$$
where $X_{n+1}=(1,x,x^2,\ldots,x^n)$, $S_{n\times (n+1)}=(a_{j-1,i-1})_{1\leq i\leq n,~1\leq j\leq n+1}$ and $a_{n,k}$ is the $(n,k)$ entry of the Riordan array $(g(t),f(t))$.\\

Also, let $(s_n(x))_{n \in \mathbb{N}}$ be the sequence of polynomials defined by equations (1.12) and (1.13), where $a_{n,k}$ is the $(n,k)$ entry of the Riordan array $(g(t),f(t))$, then

$$s_n(x)=\sum\limits_{k=0}^n b_{n,k}x^k,\eqno(1.14)$$
where $b_{n,k}$  is the $(n,k)$ entry of the Riordan array $(1/g(\bar{f}(t)),\bar{f}(t))$ and $(s_n(x))_{n \in \mathbb{N}}$ be Sheffer for $(g(t),f(t))$ \cite{Wang1}.\\

Motivated by the theory of Riordan arrays and relation between the Sheffer sequences
and Riordan array, in this paper, a new family of the 2-iterated Sheffer polynomials is introduced by means of generating function.
The determinantal definition of the 2-iterated Sheffer polynomials is established by using the relation between Sheffer sequences
and Riordan array. The quasi-monomial properties of these polynomials are derived. The 2-iterated associated Sheffer polynomials
are deduced and their properties are considered. Examples of some members belonging to these families are given.\\

\noindent
{\large{\bf 2.~2-iterated Sheffer polynomials}}
\vspace{.35cm}

The 2-iterated Sheffer polynomials (2ISP in the following) are introduced by means of generating function. Further, a
determinantal definition of the 2ISP is given.\\

In order to obtain the generating function of the 2ISP, we prove the following result:\\

\noindent
{\bf Theorem~2.1.}~{\em The 2ISP are defined by the following generating function:}
$$\frac{1}{g_1(\bar{f}_1(t))}~\frac{1}{g_2(\bar{f}_2(\bar{f}_1(t)))}\exp(x\bar{f}_2(\bar{f}_1(t)))=\sum\limits_{n=0}^\infty s_n^{[2]}(x)\frac{t^n}{n!}.\eqno(2.1)$$

\noindent
{\bf Proof.}~Let $s_n^{(1)}(x)$ and $s_n^{(2)}(x)$ be Sheffer for $(g_1(t),f_1(t))$ and $(g_2(t),f_2(t))$, respectively be two different polynomials defined by the generating functions of the forms:

$$\frac{1}{g_1(\bar{f}_1(t))}\exp(x(\bar{f}_1(t)))=\sum\limits_{n=0}^\infty s_n^{(1)}(x)\frac{t^n}{n!}\eqno(2.2a)$$
and
$$\frac{1}{g_2(\bar{f}_2(t))}\exp(x(\bar{f}_2(t)))=\sum\limits_{n=0}^\infty s_n^{(2)}(x)\frac{t^n}{n!},\eqno(2.2b)$$
respectively.\\

Expanding the exponential function and then replacing the powers $x^0,~x^1,\ldots,x^n$
by the polynomials $s_0^{(2)}(x),~s_1^{(2)}(x),\ldots,s_n^{(2)}(x)$, respectively in both sides of equation (2.2a), so that we have
$$\frac{1}{g_1(\bar{f}_1(t))}\left[1+s_1^{(2)}(x)\frac{\bar{f}_1(t)}{1!}+\ldots+s_n^{(2)}(x)\frac{(\bar{f}_1(t))^n}{n!}\right]=\sum\limits_{n=0}^\infty s_n^{(1)}(s_1^{(2)}(x))\frac{t^n}{n!}.\eqno(2.3)$$

Using equation (2.2b) with $t$ replaced by $\bar{f}_1(t)$ in the l.h.s. and denoting the resultant 2ISP in the r.h.s. of equation (2.3) by
$$s_n^{[2]}(x)=s_n^{(1)}(s_1^{(2)}(x)),\eqno(2.4)$$
we get assertion (2.1).\\

\noindent
{\bf Remark~2.1.}~We remark that equation (2.4) is the operational correspondence between the
2ISP $s_n^{[2]}(x)$ and Sheffer polynomials $s_n^{(1)}(x)$.\\

\noindent
{\bf Remark~2.2.}~We know that the Sheffer sequence for $(g(t),t)$ becomes the Appell sequence $A_n(x)$. Therefore, taking $f_1(t)=f_2(t)=t$ which gives
$\frac{1}{g_1(\bar{f}_1(t))}=\frac{1}{g_1(t)}$ and $\frac{1}{g_2(\bar{f}_2(\bar{f}_1(t)))}=\frac{1}{g_2(t)}$ in
equation (2.1) yields the generating function for the 2-iterated Appell polynomials (2IAP) $A_n^{[2]}(x)$ \cite{Iterated}.\\

\noindent
{\bf Remark~2.3.}~We know that the Sheffer sequence for $(1,f(t))$ becomes the associated Sheffer sequence $\widetilde{s}_n(x)$. Therefore, taking $g_1(t)=g_2(t)=1$ which gives
$\frac{1}{g_1(\bar{f}_1(t))}=\frac{1}{g_2(\bar{f}_2(\bar{f}_1(t)))}=1$ in
equation (2.1) yields the following consequence of Theorem 2.1:\\

\noindent
{\bf Corollary~2.1.}~{\em The 2-iterated associated Sheffer polynomials (2IASP) $\widetilde{s}_n^{[2]}(x)$ are defined by the following generating function:}
$$\exp(x\bar{f}_2(\bar{f}_1(t)))=\sum\limits_{n=0}^\infty \widetilde{s}_n^{[2]}(x)\frac{t^n}{n!}.\eqno(2.5)$$

\noindent
{\bf Remark~2.4.}~The following operational correspondence
between the 2IASP $\widetilde{s}_n^{[2]}(x)$ and associated Sheffer sequences $\widetilde{s}_n(x)$ holds:
$$\widetilde{s}_n^{[2]}(x)=\widetilde{s}_n^{(1)}(\widetilde{s}_1^{(2)}(x)).\eqno(2.6)$$

It is shown in \cite{Wang} that for $s_n^{(1)}(x)$ and $s_n^{(2)}(x)$ be Sheffer for $(g_1(t),f_1(t))$ and $(g_2(t),f_2(t))$, respectively be two different
Sheffer sequences defined by 
$$s_n^{(1)}(x)=\sum\limits_{k=0}^n b_{n,k}x^k\eqno(2.7a)$$
and
$$s_n^{(2)}(x)=\sum\limits_{k=0}^n d_{n,k}x^k,\eqno(2.7b)$$
where $b_{n,k}$ is the $(n,k)$ entry of the Riordan array $(1/g_1(\bar{f}_1(t)),\bar{f}_1(t))$ and $d_{n,k}$ is the $(n,k)$ entry of the Riordan array $(1/g_2(\bar{f}_2(t)),\bar{f}_2(t))$, then the umbral composition of $s_n^{(1)}(x)$ and $s_n^{(2)}(x)$ is the sequence $(s_n^{(2)}(s^{(1)}(x)))_{n \in \mathbb{N}}$ defined by
$$s_n^{(2)}(s^{(1)}(x))=\sum\limits_{k=0}^n d_{n,k}s_k(x)=\sum\limits_{k=0}^n d_{n,k}\sum\limits_{j=0}^k b_{k,j}x^j=\sum\limits_{j=0}^n \left(\sum\limits_{k=j}^n d_{n,k}b_{k,j}\right)x^j=\sum\limits_{j=0}^np_{n,j}x^j,\eqno(2.8)$$
where $p_{n,j}$ is the $(n,j)$ entry of the Riordan array $\left(\frac{1}{g_2(\bar{f}_2(t))g_1(\bar{f}_1(\bar{f}_2(t)))},\bar{f}_1(\bar{f}_2(t))\right)$
and $s_n^{(2)}(s^{(1)}(x))$ is Sheffer for $(g_1(t)g_2(f_1(t)),f_2(f_1(t)))$.\\

\noindent
{\bf Remark~2.5.}~We remark  that, the 2ISP $s_n^{[2]}(x)$ are actually the composition of $s_n^{(1)}(x)$ and $s_n^{(2)}(x)$ and are defined by 
$$s_n^{[2]}(x)=\sum\limits_{k=0}^n d_{n,k}s_k(x),\eqno(2.9)$$
where $d_{n,k}$ is the is the $(n,k)$ entry of the Riordan array $\left(\frac{1}{g_2(\bar{f}_2(t))g_1(\bar{f}_1(\bar{f}_2(t)))},\bar{f}_1(\bar{f}_2(t))\right)$.\\

The series definition (2.9) can also be obtained by replacing the powers $x^1$ and $x^k$
by the polynomials $s_1^{(2)}(x)$ and $s_k^{(2)}(x)$, respectively, in equation (2.7b) and then
using equation (2.4) in the l.h.s. of resultant equation.\\

Now, we derive the determinantal definition for the 2ISP $s_n^{[2]}(x)$. For this we prove the following result:\\

\noindent
{\bf Theorem~2.2.}~{\em The 2ISP $s_n^{[2]}(x)$ of degree $n$ are defined by}

$$\begin{array}{l}
s_0^{[2]}(x)=\frac{1}{a_{0,0}},\hspace{5in}(2.10)\\
s_n^{[2]}(x)=\frac{{(-1)}^n}{a_{0,0}a_{1,1}\ldots a_{n,n}}\left|\begin{array}{cccccc}
 1  & s_1^{(2)}(x)  & s_2^{(2)}(x)  & \cdots & s_{n-1}^{(2)}(x)   & s_{n}^{(2)}(x)\\
\\
 a_{0,0} &  a_{1,0}  &  a_{2,0}  & \cdots & a_{n-1,0} &a_{n,0} \\
\\
 0  &  a_{1,1}  &  a_{2,1}  & \cdots & a_{n-1,1} &a_{n,1} \\
\\
 0 & 0 & a_{2,2} & \cdots & a_{n-1,2}&  a_{n,2}\\
 . & . & . & \cdots & . & . \\
 . & . & . & \cdots & . & . \\
 0 & 0 & 0 & \cdots & a_{n-1,n-1} & a_{n,n-1}
\end{array} \right|,$$\\
\\
\hspace{.5in}=\frac{{(-1)}^n}{a_{0,0}a_{1,1}\ldots a_{n,n}} \textrm{det}\left(\begin{array}{c}
                                                                       S_{n+1}(x) \\
                                                                       M_{n \times (n+1)} \\
                                                                     \end{array}\right),\end{array}\eqno(2.11)$$
where $S_{n+1}(x)=(1,s_1^{(2)}(x),\ldots,s_{n}^{(2)}(x))$, $M_{n\times (n+1)}=(a_{j-1,i-1})_{1\leq i\leq n,~1\leq j\leq n+1}$ and $a_{n,k}$ is the $(n,k)$ entry of the Riordan array $(g_1(t)g_2(f_1(t)),f_2(f_1(t)))$.\\

\noindent
{\bf Proof.}~Replacing the power $x^n$ by the polynomial $s_n^{(2)}(x)$ in the l.h.s. and $x$ by $s_1^{(2)}(x)$ in r.h.s. of equation (1.11) and
then using operational correspondence (2.4) in the r.h.s. of resultant equation, we find
$$s_n^{(2)}(x)=\sum\limits_{k=0}^n a_{n,k} s_k^{[2]}(x),\eqno(2.12)$$
where $a_{n,k}$ is the $(n,k)$ entry of the Riordan array $(g_1(t)g_2(f_1(t)),f_2(f_1(t)))$.\\

The above identity leads to the following system of infinite equations in the unknown $s_n^{[2]}(x)$ for $n=0,1,\ldots$~:\\
$$\begin{cases}
a_{0,0}s_0^{[2]}(x)=1,\cr
\\
a_{1,0}s_0^{[2]}(x)+a_{1,1}s_1^{[2]}(x)=s_1^{(2)}(x),\cr
\\
a_{2,0}s_0^{[2]}(x)+a_{2,1}s_1^{[2]}(x)+a_{2,2}s_2^{[2]}(x)=s_2^{(2)}(x),\cr
\vdots\\
a_{n,0}s_0^{[2]}(x)+a_{n,1}s_1^{[2]}(x)+a_{n,2}s_2^{[2]}(x)+\ldots+a_{n,n}s_n^{[2]}(x)=s_n^{(2)}(x).\cr
\end{cases}\eqno(2.13)$$
\\
From first equation of system (2.13), we get assertion (2.10). Also, the special form of system (2.13) (lower triangular) allows us to work out the unknown $s_n^{[2]}(x)$. Operating with the first $n+1$ equations simply by applying the Cramer's rule, we have

$$s_n^{[2]}(x)=\frac{1}{a_{0,0}a_{1,1}\ldots a_{n,n}}\left|\begin{array}{cccccc}
                                                       a_{0,0} & 0 & 0 & \cdots & 0 & 1 \\
                                                       \\
                                                       a_{1,0}& a_{1,1} & 0 & \cdots & 0 & s_1^{(2)}(x) \\
                                                       \\
                                                       a_{2,0} & a_{2,1} & a_{2,2} & \cdots & 0 & s_2^{(2)}(x) \\
                                                       \\
                                                       . & . & . & \cdots & . & .\\
                                                       a_{n-1,0}& a_{n-1,1} & a_{n-1,2} & \cdots & a_{n-1,n-1} & s_{n-1}^{(2)}(x) \\
                                                       \\
                                                       a_{n,0} & a_{n,1} & a_{n,2} & \cdots & a_{n,n-1} &s_n^{(2)}(x)
                                                     \end{array}\right|.\eqno(2.14)$$
  \\
Then, bringing $(n+1)$-th column to the first place by $n$ transpositions of adjacent column and noting that
the determinant of a square matrix is the same as that of its transpose, we obtain assertion (2.11).\\

\noindent
{\bf Remark~2.6.}~We know that for $(1,f(t))$, the Sheffer sequences $s_n(x)$
become the associated sequences $\widetilde{s}_n(x)$, i.e., for $c_n=n!$ in equation (1.9), $(s_n(x))_{n \in \mathbb{N}}$
is associated to $f(t)$ and $a_{n,k}$  is then the $(n,k)$ entry
of the Riordan array $(1,f(t))$ and coefficients of the associated Sheffer polynomials obtained are of the form:

$$a_{n,0}=\left[\frac{t^n}{c_n}\right]\frac{(f(t))^0}{c_0}=\frac{c_n}{c_0}[t^n]1=\delta_{n,0},\eqno(2.15)$$

In view of equation (2.15), the following determinantal definition of the associated Sheffer sequences $\widetilde{s}_n(x)$ holds true \cite{Wang}:\\

$$\begin{array}{l}
\widetilde{s}_0(x)=\frac{1}{a_{0,0}},\hspace{5in}(2.16)\\
\widetilde{s}_n(x)=\frac{{(-1)}^n}{a_{0,0}a_{1,1}\ldots a_{n,n}}\left|\begin{array}{ccccc}
  x  & x^2  & \cdots & x^{n-1}   & x^n\\
\\
   a_{1,1}  &  a_{2,1}  & \cdots & a_{n-1,1} &a_{n,1} \\
\\
  0 & a_{2,2} & \cdots & a_{n-1,2}&  a_{n,2}\\
  . & . & \cdots & . & . \\
  . & . & \cdots & . & . \\
  0 & 0 & \cdots & a_{n-1,n-1} & a_{n,n-1}
\end{array} \right|,\\
\\
\hspace{.5in}=\frac{{(-1)}^n}{a_{0,0}a_{1,1}\ldots a_{n,n}} \textrm{det}\left(\begin{array}{c}
                                                                       \widetilde{S}_{n}(x) \\
                                                                       \widetilde{M}_{(n-1) \times n} \\
                                                                     \end{array}\right),\end{array}\eqno(2.17)$$
where $\widetilde{S}_{n}(x)=(x,\ldots,x^n)$, $\widetilde{M}_{(n-1)\times n}=
(a_{j,i})_{1\leq i\leq n-1,~1\leq j\leq n}$ and $a_{n,k}$ is the $(n,k)$ entry of the Riordan array $(1,f(t))$.\\

Using equation (2.15) and replacing $s_n^{(2)}(x)$ by $\widetilde{s}_n^{(2)}(x)$ in the r.h.s. of the
determinantal definition (2.10) and (2.11) of the 2ISP $s_n^{[2]}(x)$, we obtain the following consequence
of Theorem 2.2:\\

\noindent
{\bf Corollary~2.2.}~{\em The 2IASP $\widetilde{s}_n^{[2]}(x)$ of degree $n$ are defined by}

$$\begin{array}{l}
\widetilde{s}_0^{[2]}(x)=\frac{1}{a_{0,0}},\hspace{5in}(2.18)\\
\widetilde{s}_n^{[2]}(x)=\frac{{(-1)}^n}{a_{0,0}a_{1,1}\ldots a_{n,n}}\left|\begin{array}{ccccc}
  \widetilde{s}_1^{(2)}(x)  & \widetilde{s}_2^{(2)}(x)  & \cdots & \widetilde{s}_{n-1}^{(2)}(x)   & \widetilde{s}_{n}^{(2)}(x)\\
\\
   a_{1,1}  &  a_{2,1}  & \cdots & a_{n-1,1} &a_{n,1} \\
\\
  0 & a_{2,2} & \cdots & a_{n-1,2}&  a_{n,2}\\
  . & . & \cdots & . & . \\
  . & . & \cdots & . & . \\
  0 & 0 & \cdots & a_{n-1,n-1} & a_{n,n-1}
\end{array} \right|,\\
\\
\hspace{.5in}=\frac{{(-1)}^n}{a_{0,0}a_{1,1}\ldots a_{n,n}} \textrm{det}\left(\begin{array}{c}
                                                                       \widetilde{S}_{n}(x) \\
                                                                       \widetilde{M}_{(n-1) \times n} \\
                                                                     \end{array}\right),\end{array}\eqno(2.19)$$
where $\widetilde{S}_{n}(x)=(\widetilde{s}_1^{(2)}(x),\ldots,\widetilde{s}_{n}^{(2)}(x))$, $\widetilde{M}_{(n-1)\times n}=
(a_{j,i})_{1\leq i\leq n-1,~1\leq j\leq n}$ and $a_{n,k}$ is the $(n,k)$ entry of the Riordan array $(1,f_2(f_1(t)))$.\\

\noindent
{\bf Remark~2.7.}~We know that for $(g(t),t)$, the Sheffer sequences $s_n(x)$ become Appell sequences $A_n(x)$. The coefficients $a_{n,k}$ in equation (1.13) are then the $(n,k)$ entry of the exponential Riordan array $(g(t),t)$ as follows:
$$a_{n,k}=\left[\frac{t^n}{n!}\right]g(t)\frac{t^k}{k!}=\frac{n!}{k!}[t^{n-k}]g(t)={n \choose k}g_{n-k}.\eqno(2.20)$$

Therefore, using equation (2.20) in the determinantal definition (2.10) and (2.11) of the 2ISP $s_n^{[2]}(x)$ yields the determinantal
definition of the 2IAP $A_n^{[2]}(x)$, which is given in \cite{2IqAP}.\\

\noindent
{\large{\bf 3.~Quasi-monomial and other properties}}
\vspace{.35cm}

In this section, we frame the 2ISP $s_n^{[2]}(x)$ within the context of monomiality principle and
derive certain other properties of the 2ISP $s_n^{[2]}(x)$.\\

In order to derive the multiplicative and derivative operators for the 2ISP $s_n^{[2]}(x)$,
we prove the following result:\\

\noindent
{\bf Theorem~3.1.}~{\em The 2ISP $s_n^{[2]}(x)$  are quasi-monomial with respect to the following multiplicative and derivative operators:}
$$\hat{M}_{s^{[2]}}=\left(x-\frac{g'_1(f_2(\partial_x))}{g_1(f_2(\partial_x))}f'_1(f_2(\partial_x))-\frac{g'_2(\partial_x)}
{g_2(\partial_x)}\right)\frac{1}{f'_1(f_2(\partial_x))~f'_2(\partial_x)}\eqno(3.1a)$$
{\em and }
$$\hat{P}_{s^{[2]}}=f_1(f_2(\partial_x)),\eqno(3.1b)$$
{\em respectively, where $\partial_x:=\frac{\partial }{\partial x}$.}\\

\noindent
{\bf Proof.}~Differentiating equation (2.1) partially with respect to $t$, we find
\[\begin{split}
&\left(x\bar{f}'_2(\bar{f}_2(\bar{f}_1(t)))~\bar{f}'_1(\bar{f}_1(t))-\frac{g'_2(\bar{f}_2(\bar{f}_1(t)))}{g_2(\bar{f}_2(\bar{f}_1(t)))}
\bar{f}'_2(\bar{f}_2(\bar{f}_1(t)))~\bar{f}'_1(\bar{f}_1(t))-\frac{g'_1(\bar{f}_1(t))}{g_1(\bar{f}_1(t))}\bar{f}'_1(\bar{f}_1(t))\right)\\
&\frac{1}{g_1(\bar{f}_1(t))}~\frac{1}{g_2(\bar{f}_2(\bar{f}_1(t)))}\exp(x\bar{f}_2(\bar{f}_1(t)))=\sum\limits_{n=0}^\infty s_{n+1}^{[2]}(x)\frac{t^n}{n!},\hspace{2in}(3.2)\\
\end{split}\]
which can also be simplified as
\[\begin{split}
&\left(x-\frac{g'_2(\bar{f}_2(\bar{f}_1(t)))}{g_2(\bar{f}_2(\bar{f}_1(t)))}
-\frac{g'_1(\bar{f}_1(t))}{g_1(\bar{f}_1(t))}f'_2(\bar{f}_2(\bar{f}_1(t)))\right)\frac{1}{{f'}_2(\bar{f}_2({f}_1(t)))~\bar{f}'_1(\bar{f}_1(t))}
\Big\{\frac{1}{g_1(\bar{f}_1(t))}~\frac{1}{g_2(\bar{f}_2(\bar{f}_1(t)))}\\
&\exp(x\bar{f}_2(\bar{f}_1(t)))\Big\}=\sum\limits_{n=0}^\infty s_{n+1}^{[2]}(x)\frac{t^n}{n!}.\hspace{3.50in}(3.3)\\
\end{split}\]

Since $g_1(t)$ and $g_2(t)$ are invertible series of $t$, therefore $\frac{g'_1(t)}{g_1(t)}$ and $\frac{g'_2(t)}{g_2(t)}$ possess power series expansions of $t$. Thus, in view of the following identity for the 2ISP $s_n^{[2]}(x)$:

$$\partial_x\left\{\frac{1}{g_1(\bar{f}_1(t))}~\frac{1}{g_2(\bar{f}_2(\bar{f}_1(t)))}\exp(x\bar{f}_2(\bar{f}_1(t)))\right\}=\bar{f}_2(\bar{f}_1(t))~
\left(\frac{1}{g_1(\bar{f}_1(t))}~\frac{1}{g_2(\bar{f}_2(\bar{f}_1(t)))}\exp(x\bar{f}_2(\bar{f}_1(t)))\right),\eqno(3.4a)$$
or, equivalently
$$f_1(f_2(\partial_x))\left\{\frac{1}{g_1(\bar{f}_1(t))}~\frac{1}{g_2(\bar{f}_2(\bar{f}_1(t)))}\exp(x\bar{f}_2(\bar{f}_1(t)))\right\}=t~
\left(\frac{1}{g_1(\bar{f}_1(t))}~\frac{1}{g_2(\bar{f}_2(\bar{f}_1(t)))}\exp(x\bar{f}_2(\bar{f}_1(t)))\right),\eqno(3.4b)$$
equation (3.3) becomes
\[\begin{split}
&\left(x-\frac{g'_1(f_2(\partial_x))}{g_1(f_2(\partial_x))}f'_1(f_2(\partial_x))-\frac{g'_2(\partial_x)}
{g_2(\partial_x)}\right)\frac{1}{f'_1(f_2(\partial_x))~f'_2(\partial_x)}\Big\{\frac{1}{g_1(\bar{f}_1(t))}~\frac{1}{g_2(\bar{f}_2(\bar{f}_1(t)))}\\
&\exp(x\bar{f}_2(\bar{f}_1(t)))\Big\}=\sum\limits_{n=0}^\infty s_{n+1}^{[2]}(x)\frac{t^n}{n!}.\hspace{3.50in}(3.5)\\
\end{split}\]
Again, using generating function (2.1) in the l.h.s. of equation (3.5) and then rearranging the summation yields

$$
\sum\limits_{n=0}^\infty\left(\left(x-\frac{g'_1(f_2(\partial_x))}{g_1(f_2(\partial_x))}f'_1(f_2(\partial_x))-\frac{g'_2(\partial_x)}
{g_2(\partial_x)}\right)\frac{1}{f'_1(f_2(\partial_x))~f'_2(\partial_x)}\right)\left\{s_{n}^{[2]}(x)\frac{t^n}{n!}\right\}= \sum\limits_{n=0}^\infty s_{n+1}^{[2]}(x)\frac{t^n}{n!}.\eqno(3.6)$$
Equating the coefficients of the same powers of $t$ in both sides of the above equation, we find
$$\left(\left(x-\frac{g'_1(f_2(\partial_x))}{g_1(f_2(\partial_x))}f'_1(f_2(\partial_x))-\frac{g'_2(\partial_x)}
{g_2(\partial_x)}\right)\frac{1}{f'_1(f_2(\partial_x))~f'_2(\partial_x)}\right)\left\{s_{n}^{[2]}(x)\right\}= s_{n+1}^{[2]}(x),\eqno(3.7)$$
which in view of monomiality principle equation $\hat{M}\{s_n(x)\}=s_{n+1}(x)$ \cite{HerSheff} for $s_n^{[2]}(x)$ yields assertion (3.1a).\\

In order to prove assertion (3.1b), we use generating function (2.1) in both sides of the identity (3.4b), so that we have
$$f_1(f_2(\partial_x))\left\{ \sum\limits_{n=0}^\infty s_{n}^{[2]}(x)\frac{t^n}{n!}\right\}=\sum\limits_{n=1}^\infty s_{n-1}^{[2]}(x)\frac{t^n}{(n-1)!}.\eqno(3.8)$$

Rearranging the summation in the l.h.s. of equation (3.8) and then equating the coefficients of the same powers of $t$ in both sides of the resultant equation, we find
$$f_1(f_2(\partial_x))\left\{s_{n}^{[2]}(x)\right\}=n~s_{n-1}^{[2]}(x),~n\geq1,\eqno(3.9)$$
which in view of monomiality principle equation  $\hat{P}\{s_n(x)\}=n~s_{n-1}(x)$ \cite{HerSheff} for $s_n^{[2]}(x)$ yields assertion (3.1b).\\

\noindent
{\bf Theorem~3.2.}~{\em The 2ISP $s_n^{[2]}(x)$ satisfy the following differential equation:}
$$\left(\left(\frac{x}{f'_1(f_2(\partial_x))}-\frac{g'_1(f_2(\partial_x))}{g_1(f_2(\partial_x))}\right)\frac{f_1(f_2(\partial_x))}{f'_2(\partial_x)}-\frac{g'_2(\partial_x)}
{g_2(\partial_x)}\frac{f_1(f_2(\partial_x))}{f'_1(f_2(\partial_x))~f'_2(\partial_x)}-n\right)s_n^{[2]}(x)=0.\eqno(3.10)$$

\noindent
{\bf Proof}.~Using equations (3.1a) and (3.1b) in monomiality principle equation $\hat{M}\hat{P}\{s_n(x)\}=ns_{n}(x)$ \cite{HerSheff} for $s_n^{[2]}(x)$, we get assertion (3.10).\\

\noindent
{\bf Remark~3.1.}~We know that the Sheffer sequence for $(1,f(t))$ becomes the associated Sheffer sequence $\widetilde{s}_n(x)$. Therefore, taking $g_1(t)=g_2(t)=1$ which implies $g'_1(t)=g'_2(t)=0$ in equations (3.1a), (3.1b) and (3.10), we get the following consequence of Theorem 3.1 and  Theorem 3.2:\\

\noindent
{\bf Corollary~3.1.}~{\em The 2IASP $\widetilde{s}_n^{[2]}(x)$  are quasi-monomial with respect to the following multiplicative and derivative operators:}
$$\hat{M}_{\widetilde{s}^{[2]}}=\frac{x}{{f'_1(f_2(\partial_x))~f'_2(\partial_x)}}\eqno(3.11a)$$
{\em and }
$$\hat{P}_{\widetilde{s}^{[2]}}=f_1(f_2(\partial_x)),\eqno(3.11b)$$
{\em respectively.}\\

\noindent
{\bf Corollary~3.2.}~{\em The 2IASP $\widetilde{s}_n^{[2]}(x)$ satisfy the following differential equation:}
$$\left(\frac{x~f_1(f_2(\partial_x))}{{f'_1(f_2(\partial_x))~f'_2(\partial_x)}}-n\right)\widetilde{s}_n^{[2]}(x)=0.\eqno(3.12)$$

\noindent
{\bf Remark~3.2.}~We know that for $(g(t),t)$, the Sheffer sequences $s_n(x)$ reduce to the Appell sequences $A_n(x)$.
Therefore, taking $f_1(t)=f_2(t)=t$ which implies $f'_1(t)=f'_2(t)=1$ in equations (3.1a), (3.1b) and (3.10) yields the multiplicative and derivative
operators and differential equation for the 2IAP $A_n^{[2]}(x)$, which are given in \cite{Iterated}.\\

Now, we prove the conjugate representation of the 2ISP $s_n^{[2]}(x)$. For this, we prove the following result:\\

\noindent
{\bf Theorem~3.2.}~{\em Let $s_n^{[2]}(x)$ be Sheffer for $(g_1(t)g_2(f_1(t)),f_2(f_1(t)))$, then the following
conjugate representation for the 2ISP $s_n^{[2]}(x)$ holds true:}

$$s_n^{[2]}(x)=\sum\limits_{k=0}^n \frac{\langle(g_2(\bar{f}_2(t))g_1(\bar{f}_1(\bar{f}_2(t))))^{-1}(\bar{f}_1(\bar{f}_2(t)))^k \mid s_n(x)\rangle}{c_k}s_k(x).\eqno(3.13)$$

\noindent
{\bf Proof.}~Since $s_n^{[2]}(x)$ are defined by equation (2.9), where $d_{n,k}$ is the $(n,k)$ entry of the Riordan array
$\left(\frac{1}{g_2(\bar{f}_2(t))g_1(\bar{f}_1(\bar{f}_2(t)))},\bar{f}_1(\bar{f}_2(t))\right)$. Then, according to the
definition of Riordan arrays, we have
$$d_{n,k}=\left[\frac{t^n}{c_n}\right]\frac{1}{g_2(\bar{f}_2(t))g_1(\bar{f}_1(\bar{f}_2(t)))}\frac{(\bar{f}_1(\bar{f}_2(t)))^k}
{c_k}.\eqno(3.14)$$

Simplifying the above equation yields
\[\begin{split}
d_{n,k}&=\frac{1}{c_k}\left[\frac{t^n}{c_n}\right](g_2(\bar{f}_2(t))g_1(\bar{f}_1(\bar{f}_2(t))))^{-1}((\bar{f}_1(\bar{f}_2(t))))^k\\
&=\frac{1}{c_k}\langle(g_2(\bar{f}_2(t))g_1(\bar{f}_1(\bar{f}_2(t))))^{-1}(\bar{f}_1(\bar{f}_2(t)))^k \mid s_n(x)\rangle,\hspace{2in}(3.15)\\
\end{split}\]

Now, using equation (3.15) in equation (2.9), we are led to assertion (3.13).\\

In the next section, examples of some members belonging to the 2ISP and 2IASP are considered.\\

\noindent
{\large{\bf 4.~Examples}}
\vspace{.35cm}

We recall that the Laguerre polynomials $L_n^{(\alpha)}(x)$
of order $\alpha$ form the Sheffer sequence for the pair $\left(\frac{1}{(1-t)^{\alpha+1}},\frac{t}{t-1}\right)$ \cite{Andrew,Roman,Rainville} are defined by the generating function:
$$\frac{1}{(1-t)^{\alpha+1}}\exp\Big(\frac{-xt}{1-t}\Big)=\sum\limits_{n=0}^\infty L_n^{(\alpha)}(x)t^n,\eqno(4.1)$$
which for $\alpha=0$ gives the generating function of the Laguerre polynomials $L_n(x)$ as \cite{Andrew}:
$$\frac{1}{(1-t)}\exp\Big(\frac{-xt}{1-t}\Big)=\sum\limits_{n=0}^\infty L_n^{(\alpha)}(x)t^n.\eqno(4.2)$$

Also, the polynomials of the Gegenbauer case denoted by ${-\lambda \choose n}s_n(x)$ considered in \cite{Wang} form the Sheffer sequence for the pair
$\left(\left(\frac{2}{1+\sqrt{1-t^2}}\right)^{\lambda_0},\frac{-t}{1+\sqrt{1-t^2}}\right)$ and are defined by the
generating function of the form:
$$(1+t^2)^{\lambda-\lambda_0}(1-2xt+t^2)^{-\lambda}=\sum\limits_{n=0}^\infty {-\lambda \choose n}s_n(x)t^n,\eqno(4.3)$$
which for $\lambda_0=\lambda$ gives the generating function of the Gegenbauer polynomials $C_n^{(\lambda)}(x)$ as \cite{Abram,Rainville,Roman1}:
$$(1-2xt+t^2)^{-\lambda}=\sum\limits_{n=0}^\infty C_n^{(\lambda)}(x)t^n.\eqno(4.4)$$

We note that corresponding to each member belonging to the Sheffer family, there exists a new special polynomial belonging to the 2ISP family. Thus, by making
suitable choice for the functions $g(t)$ and $f(t)$ in equation (2.1), we get the generating function for the corresponding member belonging to the 2ISP
family. The other properties of these special polynomials can be obtained from the results derived in previous sections.\\

We consider the following examples:\\

\noindent
{\bf Example~4.1.}~Taking $g_1(t)=g_2(t)=\frac{1}{(1-t)^{\alpha+1}}$ and $f_1(t)=f_2(t)=\frac{t}{t-1}$ in the l.h.s. of generating function (2.1), we find that the resultant 2-iterated Laguerre polynomials (2ILP) of order $\alpha$, which may be denoted by $L_n^{(\alpha)[2]}(x)$ in the r.h.s. are defined by the following generating function:
$$\frac{1}{(1-t)^{\alpha+1}(1-t)^{-\alpha-1}}\exp(xt)=\sum\limits_{n=0}^\infty L_n^{(\alpha)[2]}(x){t^n}.\eqno(4.5)$$

The operational correspondence between the 2ILP of order $\alpha$ $L_n^{(\alpha)[2]}(x)$ and Laguerre polynomials of order, $\alpha$ $L_n^{(\alpha)}(x)$ is given by:
$$L^{(\alpha)[2]}(x)=L_n^{(\alpha,1)}(L_1^{(\alpha,2)}(x)).\eqno(4.6)$$

The 2ILP of order $\alpha$, $L_n^{(\alpha)[2]}(x)$ are quasi-monomial with respect to the following multiplicative and
derivative operators:
$$\hat{M}_{L^{(\alpha)[2]}}=\left(x-(\alpha+1)(\partial_x -1)^{3}-\frac{(\alpha+1)}{(\partial_x -1)}\right),\eqno(4.7a)$$
$$\hat{P}_{L^{(\alpha)[2]}}=\partial_x.\eqno(4.7b)$$

The 2ILP of order $\alpha$, $L_n^{(\alpha)[2]}(x)$ satisfy the following differential equation:
$$\left(x\partial_x-(\alpha+1)(\partial_x -1)^{3}\partial_x-\frac{(\alpha+1)\partial_x}{(\partial_x -1)}-n\right)L_n^{(\alpha)[2]}(x)=0.\eqno(4.8)$$

Consider the following series definition for the Laguerre polynomials of order $\alpha$, $L_n^{(\alpha)}(x)$ \cite{Roman}:
$$L_n^{(\alpha)}(x)=\sum\limits_{k=0}^n  \frac{(-1)^k}{k!}{n+\alpha \choose n-k}x^k.\eqno(4.9)$$

Replacing the power $x^k$ by the polynomial $L_k^{(\alpha)}(x)$ in the r.h.s. and $x$ by $L_1^{(\alpha)}(x)$ in the l.h.s. of equation
(4.9) and then using equation (4.6) in the l.h.s. of resultant equation yields the following series definition for the 2ILP of order $\alpha$, $L_n^{(\alpha)[2]}(x)$:
$$L_n^{(\alpha)[2]}(x)=\sum\limits_{k=0}^n \frac{(-1)^k }{k!}{n+\alpha \choose n-k}L_k^{(\alpha)}(x).\eqno(4.10)$$

It has been shown in \cite{Wang} that for $a_{n,k}=(-1)^k \frac{n!}{k!}{n+\alpha \choose n-k}$ the determinantal definition of the
Sheffer polynomials given by equations (1.12) and (1.13) reduces to the determinantal definition of the Laguerre polynomials of order $\alpha$, $L_n^{(\alpha)}(x)$.\\

Therefore, taking $s_n(x)=L_n^{(\alpha)}(x)$ and $a_{n,k}=(-1)^k \frac{n!}{k!}{n+\alpha \choose n-k}$ in equations (2.10) and (2.11), we find that the following determinantal definition of the 2ILP of order $\alpha$, $L_n^{(\alpha)[2]}(x)$ holds true:\\
\\
$$\begin{array}{l}
L_0^{(\alpha)[2]}(x)=1,\hspace{5in}(4.11)\\
\\

L_n^{(\alpha)[2]}(x)=(-1)^{\frac{n(n+3)}{2}}
\left|\begin{array}{cccccc}
  1 & L_1^{(\alpha)}(x) & L_2^{(\alpha)}(x) & . & L_{n-1}^{(\alpha)}(x) & L_{n}^{(\alpha)}(x) \\
  \\
  1 & \alpha+1 & (\alpha+2)_2 &  . & (\alpha+n-1)_{n-1} & (\alpha+n)_n  \\
  \\
  0 & -1 & -2(\alpha+2) & . & -(n-1)(\alpha+n-1)_{n-2} & -n(\alpha+n)_{n-1} \\
  \\
  0 & 0 & 1 & .  & \frac{(n-1)(n-2)}{2}(\alpha+n-1)_{n-3} & \frac{n(n-1)}{2}(\alpha+n)_{n-2} \\

  \vdots & \vdots & \vdots & \vdots  & \vdots & \vdots \\
 \\
  0 & 0 & 0 & . & (-1)^{n-1} & (-1)^{n-1}n(n+\alpha)  \\
 \end{array}\right|,\\
\\
\hspace{.5in}= (-1)^{\frac{n(n+3)}{2}}\textrm{det}\left(\begin{array}{c}
                                                                       S_{n+1}(x) \\
                                                                       M_{n \times (n+1)} \\
                                                                     \end{array}\right),\end{array}\eqno(4.12)$$
where $S_{n+1}(x)=(1,L_1^{(\alpha)}(x),\ldots,L_{n}^{(\alpha)}(x))$, $M_{n\times (n+1)}=(a_{j-1,i-1})_{1\leq i\leq n,~1\leq j\leq n+1}$
and $L_n^{(\alpha)}(x)~(n=0,1,\ldots)$ are the Laguerre polynomials of order $\alpha$.\\

\noindent
{\bf Remark~4.1.}~We remark that by taking $\alpha=0$ in the results derived above for the 2-iterated Laguerre polynomials (2ILP) of order $\alpha$,
we get the corresponding results for the  2-iterated Laguerre polynomials (2ILP), which may be denoted by $L_{n}^{[2]}(x)$.\\

\noindent
{\bf Example~4.2.}~Taking $g_1(t)=g_2(t)=\left(\frac{2}{1+\sqrt{1-t^2}}\right)^{\lambda_0}$ and $f_1(t)=f_2(t)=\frac{-t}{1+\sqrt{1-t^2}}$ in the l.h.s. of generating function (2.1), we find that the resultant 2-iterated polynomials of the Gegenbauer case (2IPoGc), which may be denoted by ${-\lambda \choose n}\mathfrak{s}_n^{[2]}(x)$ in the r.h.s. are defined by the following generating function:

$$\left(\frac{1+t^2}{1+6t^2+t^4}\right)^{\lambda_0}\exp{\left(\frac{4xt(1+t^2)}{1+6t^2+t^4}\right)}=\sum\limits_{n=0}^\infty {-\lambda \choose n}\mathfrak{s}_n^{[2]}(x)t^n.\eqno(4.13)$$

The operational correspondence between the 2IPoGc $\mathfrak{s}_n^{[2]}(x)$ and polynomials of the Gegenbauer case $\mathfrak{s}_n(x)$
is given by:
$$\mathfrak{s}_n^{[2]}(x)=\mathfrak{s}_n^{(1)}(\mathfrak{s}_1^{(2)}(x)).\eqno(4.14)$$

The 2IPoGc ${-\lambda \choose n}\mathfrak{s}_n^{[2]}(x)$ are quasi-monomial with respect to the following multiplicative and derivative operators:
\[\begin{split}
&\hat{M}_{\mathfrak{s}^{(\alpha)[2]}}=\Big(x-\frac{\lambda_0 2^{\lambda_0 -1}\partial_x (1+\sqrt{1-\partial_x^2})^{\lambda_{0} +2}}{\big(\sqrt{2(1-\partial_x^2+\sqrt{1-\partial_x^2})}\big)^2~\big(1+\sqrt{1-\partial_x^2}+\sqrt{2(1-\partial_x^2+\sqrt{1-\partial_x^2})\big)^{\lambda_{0}+1}}}\\
&-\frac{\lambda_0 \partial_x}{\sqrt{1-\partial_x^2}(1+\sqrt{1-\partial_x^2})}\Big)\frac{1+\sqrt{1-\partial_x^2}}
{\sqrt{2(1-\partial_x^2+\sqrt{1-\partial_x^2})}~(1+\sqrt{1-\partial_x^2}+\sqrt{2(1-\partial_x^2+\sqrt{1-\partial_x^2})})}.\\
&\times \frac{1}{\sqrt{(1-\partial_x^2)}}\hspace{5in}(4.15)\\
\end{split}\]

$$\hat{P}_{\mathfrak{s}^{(\alpha)[2]}}=\frac{-\partial_x}{1+\sqrt{1-\partial_x^2}+\sqrt{2(1-\partial_x^2+\sqrt{1-\partial_x^2})}}.\eqno(4.16)$$

The 2IPoGc ${-\lambda \choose n}\mathfrak{s}_n^{[2]}(x)$ satisfy the following differential equation:
\[\begin{split}
&\Big(\Big(x+\frac{\lambda_0 2^{\lambda_0 -1}\partial_x (1+\sqrt{1-\partial_x^2})^{\lambda_{0} +2}}{\big(\sqrt{2(1-\partial_x^2+\sqrt{1-\partial_x^2})}\big)^2~\big(1+\sqrt{1-\partial_x^2}+\sqrt{2(1-\partial_x^2+\sqrt{1-\partial_x^2})\big)^{\lambda_{0}+1}}}\\
&+\frac{\lambda_0 \partial_x}{\sqrt{1-\partial_x^2}(1+\sqrt{1-\partial_x^2})}\Big)\frac{\partial_x(1+\sqrt{1-\partial_x^2})}
{\sqrt{2(1-\partial_x^2+\sqrt{1-\partial_x^2})}~(1+\sqrt{1-\partial_x^2}+\sqrt{2(1-\partial_x^2+\sqrt{1-\partial_x^2})})^2}\\
&\times \frac{1}{\sqrt{1-\partial_x^2}}-n\Big)\mathfrak{s}_n^{[2]}(x)=0.\hspace{4in}(4.17)\\
\end{split}\]
\\
It has been shown in \cite{Wang} that by taking
$$a_{n,k}=\begin{cases}
0\hspace{3cm}\textrm{if}~n-k~\textrm{is odd},\cr
\frac{c_n (-1)^k \lambda_0 +k}{c_k 2^n \lambda_0 +n}{\lambda_0 +n \choose \frac{n-k}{2}}~~\textrm{if}~n-k~\textrm{is even},\cr
\end{cases}\eqno(4.18)$$
where $c_n=1/{\lambda \choose n}$, the determinantal definition of the
Sheffer polynomials given by equations (1.12) and (1.13) reduces to the determinantal definition of the polynomials of Gegenbauer case ${-\lambda \choose n}\mathfrak{s}_n(x)$.\\

Therefore, taking $s_n^{(2)}(x)=\mathfrak{s}_n{x}$ and using equation (4.18) in equations (2.10) and (2.11), we find that the following determinantal definition of the 2IPoGc ${-\lambda \choose n}\mathfrak{s}_n^{[2]}(x)$ holds true:\\

$$\begin{array}{l}
\mathfrak{s}_0^{[2]}(x)=1,\hspace{5in}(4.19)\\
\\

\mathfrak{s}_n^{[2]}(x)=(-1)^{\frac{n(n+3)}{2}}2^{\frac{n(n+1)}{2}}
\left|\begin{array}{cccccc}
  1 & \mathfrak{s}_1(x) & \mathfrak{s}_2(x) & . & \mathfrak{s}_{n-1}(x) & \mathfrak{s}_{n}(x) \\
  \\
  1 & 0 & \frac{\lambda_0}{2\lambda(\lambda+1)} &  . & 0 & \frac{n! \lambda_0 (\lambda_0 +n-1)!}{2^n \big(\frac{n}{2}\big)!\lambda (\lambda+1)\ldots(\lambda+n-1)\big(\lambda_0 +\frac{n}{2}\big)!}\\
  \\
  0 & -\frac{1}{2} & 0 & . &  & 0 \\
  \\
  0 & 0 & \frac{1}{4} & .  & 0 & \frac{n(n-1)\ldots3(\lambda_0 +2)(\lambda_0 +n-1)!}{2^n \big(\frac{n-2}{2}\big)!(\lambda+2)\ldots(\lambda+n-1)\big(\lambda_0 +\frac{n+2}{2}\big)!}\\

  \vdots & \vdots & \vdots & \vdots  & \vdots & \vdots \\
 \\
  0 & 0 & 0 & . & \big(\frac{-1}{2}\big)^{n-1} & 0  \\
 \end{array}\right|,\\
\\
\hspace{.5in}=(-1)^{\frac{n(n+3)}{2}}2^{\frac{n(n+1)}{2}}\textrm{det}\left(\begin{array}{c}
                                                                       S_{n+1}(x) \\
                                                                       M_{n \times (n+1)} \\
                                                                     \end{array}\right),\end{array}\eqno(4.20)$$
where $S_{n+1}(x)=(1,\mathfrak{s}_1(x),\ldots,\mathfrak{s}_{n}(x))$, $M_{n\times (n+1)}=(a_{j-1,i-1})_{1\leq i\leq n,~1\leq j\leq n+1}$
and $\mathfrak{s}_n(x)~(n=0,1,\ldots)$ are the polynomials of Gegenbauer case.\\

\noindent
{\bf Remark~4.2.}~We remark that, by taking $\lambda_0=\lambda$ in the results derived above for the 2IPoGc ${-\lambda \choose n}\mathfrak{s}_n^{[2]}(x)$
we get the corresponding results for the  2-iterated Gegenbauer polynomials (2IGP) which may be denoted by $C_n^{(\lambda)[2]}(x)$.\\

\noindent
{\bf Remark~4.3.}~It is given in \cite{Abram} that, for $\lambda=1$, the Gegenbauer polynomials becomes the Chebyshev polynomials of the second kind $U_n(x)$,
i.e., $U_n(x)=C_n^{(1)}(x)$ and for $\lambda=1/2$, the Gegenbauer polynomials becomes the Legendre polynomials $P_n(x)$,
i.e., $P_n(x)=C_n^{(1/2)}(x)$. Therefore, taking $\lambda=1$ and $\lambda=1/2$ in the results of the 2-iterated Gegenbauer polynomials (2IGP) $C_n^{(\lambda)[2]}(x)$, we obtained the corresponding results of the 2-iterated Chebyshev polynomials of the second kind $U_n^{[2]}(x)$ and
2-iterated Legendre polynomials $P_n^{[2]}(x)$.\\

Now, we proceed to introduce certain members belonging to the 2IASP. We recall that the associated Sheffer family contains the falling factorials $\Big(\frac{x}{a}\Big)_n$ and exponential polynomials $\phi_n(x)$ as the important members.\\

The falling factorial $\Big(\frac{x}{a}\Big)_n$ associated to $f(t)=e^{at}-1$ \cite{Roman} is defined by the generating function:
$$\exp(xa^{-1}\log(1+t))=\sum\limits_{n=0}^\infty \left(\frac{x}{a}\right)_n\frac{t^n}{n!}\eqno(4.21)$$
and the exponential polynomials $\phi_n(x)$ associated to $f(t)=\log(1+t)$ \cite{Bell,Roman} defined by the generating function:
$$\exp(x(e^{t}-1))=\sum\limits_{n=0}^\infty \phi_n(x)\frac{t^n}{n!}.\eqno(4.22)$$

Also, for each member belonging to the associated Sheffer family, there exists a new special polynomial belonging to the 2IASP family. Thus, by making
suitable choice for the functions $f_1(t)$ and $f_2(t)$ in equation (2.5), we get the generating function for the corresponding member belonging to the 2IASP
family. The other properties of these special polynomials can be obtained from the results derived in previous sections.\\

We consider the following examples:\\

\noindent
{\bf Example~4.3.}~Taking $f_1(t)=f_2(t)=e^{at}-1$ in the l.h.s. of generating function (2.5), we find that the resultant 2-iterated falling factorial (2IFF), denoted by $\big(\frac{x}{a}\big)_n^{[2]}$ in the r.h.s. are defined by the following generating function:
$$\exp(xa^{-1}\log(1+a^{-1}\log(1+t)))=\sum\limits_{n=0}^\infty \Big(\frac{x}{a}\Big)_n^{[2]} \frac{t^n}{n!}.\eqno(4.23)$$

The operational correspondence between the 2IFF $\Big(\frac{x}{a}\Big)_n^{[2]}$ and falling factorial $\Big(\frac{x}{a}\Big)_n$ is given by:
$$\Big(\frac{x}{a}\Big)_n^{[2]}=\Big(\frac{\big(\frac{x}{a}\big)_1^{(2)}}{a}\Big)_n^{(1)}.\eqno(4.24)$$

The 2IFF $\Big(\frac{x}{a}\Big)_n^{[2]}$ are quasi-monomial with respect to the following multiplicative and derivative operators:
$$\hat{M}_{\big(\frac{x}{a}\big)_n^{[2]}}=\frac{x}{ae^{ae^{a\partial_x}}ae^{a\partial_x}},\eqno(4.25a)$$
$$\hat{P}_{\big(\frac{x}{a}\big)_n^{[2]}}=e^{a(e^{a\partial_x}-1)}-1.\eqno(4.25b)$$

The 2IFF $\big(\frac{x}{a}\big)_n^{[2]}$ satisfy the following differential equation:
$$\left(\frac{x~e^{a(e^{a\partial_x}-1)}-1}{ae^{ae^{a\partial_x}}ae^{a\partial_x}}-n\right)\Big(\frac{x}{a}\Big)_n^{[2]}=0.\eqno(4.26)$$

Consider the following series definition for the falling factorial $\Big(\frac{x}{a}\Big)_n$ \cite{Roman}:
$$\Big(\frac{x}{a}\Big)_n=\sum\limits_{k=0}^n a^n s(n,k)x^k,\eqno(4.27)$$
where $s(n,k)$ are the Stirling numbers of the first kind defined by $s(n,k)=\left[
                                                                               \begin{array}{c}
                                                                                 n \\
                                                                                 k \\
                                                                               \end{array}
                                                                             \right],~k,~n \in N,~1\leq k<n$.\\

Replacing the power $x^k$ by the polynomial $\Big(\frac{x}{a}\Big)_k$ in the r.h.s. and $x$ by $\Big(\frac{x}{a}\Big)_1$ in the l.h.s. of equation
(4.27) and then using equation (4.24) in the l.h.s. of resultant equation yields the following series definition for the 2IFF $\big(\frac{x}{a}\big)_n^{[2]}$:
$$\Big(\frac{x}{a}\Big)_n^{[2]}=\sum\limits_{k=0}^n a^n s(n,k)\Big(\frac{x}{a}\Big)_k.\eqno(4.28)$$

It has been shown in \cite{Wang}, that for $a_{n,k}=a^n S(n,k)$, where $S(n,k)$ are the Stirling numbers
of the second kind in equations (2.16) and (2.17), the determinantal definition of the associated Sheffer polynomials reduces to the
determinantal definition of falling factorials $\Big(\frac{x}{a}\Big)_n$.\\

Therefore, taking $\widetilde{s}_n(x)=\big(\frac{x}{a}\big)_n$ and $a_{n,k}=a^n S(n,k)$ in the r.h.s. of equations (2.18) and (2.19), we find that the
2IFF $\Big(\frac{x}{a}\Big)_n^{[2]}$ are defined by the following determinantal definition:

$$\begin{array}{l}
\Big(\frac{x}{a}\Big)_0^{[2]}=1,\hspace{5in}(4.29)\\
\Big(\frac{x}{a}\Big)_n^{[2]}=\frac{{(-1)}^{n+1}}{a^{{n+1 \choose 2}}}\left|\begin{array}{ccccc}
  \Big(\frac{x}{a}\Big)_1  &  \Big(\frac{x}{a}\Big)_2 & \cdots & \Big(\frac{x}{a}\Big)_{n-1}   & \Big(\frac{x}{a}\Big)_{n}\\
\\
   a S(1,1)  &  a^2 S(2,1)  & \cdots & a^{n-1} S(n-1,1) & a^{n} S(n,1) \\
\\
  0 & a^2 S(2,2) & \cdots & a^{n-1} S(n-1,2)&   a^{n} S(n,2) \\
  . & . & \cdots & . & . \\
  . & . & \cdots & . & . \\
  0 & 0 & \cdots &a^{n-1} S(n-1,n-1)  & a^{n} S(n,n-1)
\end{array} \right|,\\
\\
\hspace{.5in}=\frac{{(-1)}^{n+1}}{a^{{n+1 \choose 2}}} \textrm{det}\left(\begin{array}{c}
                                                                       \widetilde{S}_{n}(x) \\
                                                                       \widetilde{M}_{(n-1) \times n} \\
                                                                     \end{array}\right),\end{array}\eqno(4.30)$$
where $\widetilde{S}_{n}(x)=\Big(\Big(\frac{x}{a}\Big)_1,\ldots,\Big(\frac{x}{a}\Big)_{n}\Big)$, $\widetilde{M}_{(n-1)\times n}=
(a_{j,i})_{1\leq i\leq n-1,~1\leq j\leq n}$ and $\Big(\frac{x}{a}\Big)_n~(n=0,1,\ldots)$ are the falling factorials.\\

\noindent
{\bf Example~4.4.}~Taking $f_1(t)=f_2(t)=\log(1+t)$ in the l.h.s. of generating function (2.5), we find that the resultant 2-iterated exponential
polynomials (2IEP) denoted by $\phi_n^{[2]}(x)$ in the r.h.s. are defined by the following generating function:
$$e^{(e^{t}-1)}-1=\sum\limits_{n=0}^\infty \phi_n^{[2]}(x)\frac{t^n}{n!}.\eqno(4.31)$$

The operational correspondence between the 2IEP $\phi_n^{[2]}(x)$ and exponential polynomials $\phi_n(x)$ is given by:
$$\phi_n^{[2]}(x)=\phi_n^{(1)}(\phi_1^{(2)}(x)).\eqno(4.32)$$

The 2IEP $\phi_n^{[2]}(x)$ are quasi-monomial with respect to the following multiplicative and derivative operators:
$$\hat{M}_{\phi^{[2]}}=\frac{x}{(1+\log(1+t))(1+t)},\eqno(4.33a)$$
$$\hat{P}_{\phi^{[2]}}=\log(1+\log(1+t)).\eqno(4.33b)$$

The 2IEP $\phi_n^{[2]}(x)$ satisfy the following differential equation:
$$\left(\frac{x \log(1+\log(1+t))}{(1+\log(1+t))(1+t)}-n\right)\phi_n^{[2]}(x)=0.\eqno(4.34)$$

Consider the following series definition for the exponential polynomials $\phi_n(x)$ \cite{Roman}:
$$\phi_n(x)=\sum\limits_{k=0}^n S(n,k)x^k,\eqno(4.35)$$
where $S(n,k)$ are the Stirling numbers of the second kind defined by
$$ S(n,k)=\left\{
           \begin{array}{c}
             n \\
             k \\
           \end{array}
         \right\}=\frac{1}{k!}\sum\limits_{j=0}^k (-1)^{k-j}{k \choose j}j^n.\eqno(4.36)$$

Replacing the power $x^k$ by the polynomial $\phi_k(x)$  in the r.h.s. and $x$ by $\phi_1(x)$ in the l.h.s. of equation
(4.30) and then using equation (4.35) in the l.h.s. of resultant equation yields the following series definition for the
2IEP $\phi_n^{[2]}(x)$:
$$\phi_n^{[2]}(x)=\sum\limits_{k=0}^n S(n,k)\phi_k(x),\eqno(4.37)$$

It has been shown in \cite{Wang}, that for $a_{n,k}= s(n,k)$, where $s(n,k)$ are the Stirling numbers
of the first kind in equations (2.16) and (2.17), the determinantal definition of the associated Sheffer polynomials reduces to the
determinantal definition of exponential polynomials $\phi_n(x)$.\\

Therefore, taking $\widetilde{s}_n(x)=\phi_n(x)$ and $a_{n,k}=s(n,k)$ in equations (2.18) and (2.19), we find that the
2-iterated exponential polynomials (2IEP) $\phi_n^{[2]}(x)$ are defined by the following determinantal definition:

$$\begin{array}{l}
\phi_0^{[2]}(x)=1,\hspace{5in}(4.37)\\
\phi_n^{[2]}(x)={{(-1)}^{n+1}}\left|\begin{array}{ccccc}
 \phi_1(x)  &  \phi_2(x)  & \cdots & \phi_{n-1}(x)   & \phi_{n}(x)\\
\\
   s(1,1)  &  s(2,1)  & \cdots & s(n-1,1) & s(n,1) \\
\\
  0 & s(2,2) & \cdots & s(n-1,2)&   s(n,2) \\
  . & . & \cdots & . & . \\
  . & . & \cdots & . & . \\
  0 & 0 & \cdots &s(n-1,n-1)  & s(n,n-1)
\end{array} \right|,\\
\\
\hspace{.5in}={(-1)}^{n+1} \textrm{det}\left(\begin{array}{c}
                                                                       \widetilde{S}_{n}(x) \\
                                                                       \widetilde{M}_{(n-1) \times n} \\
                                                                     \end{array}\right),\end{array}\eqno(4.38)$$
where $\widetilde{S}_{n}(x)=\Big(\phi_1(x),\ldots,\phi_{n}(x)\Big)$, $\widetilde{M}_{(n-1)\times n}=
(a_{j,i})_{1\leq i\leq n-1,~1\leq j\leq n}$ and $\phi_n(x)~(n=0,1,\ldots)$ are the exponential polynomials.\\

Similar results can also be proved for certain other known members belonging to the Sheffer family. These members are listed in Table 1:\\

\noindent
\textbf{Table 1.~Some known Sheffer polynomials.}\\
\\
{\tiny{
\begin{tabular}{|l|l|l|l|l|}
\hline
&&&&\\
S.No. & ${A}(t); H(t)$ & ${g}(t); f(t)$ & Generating Functions & Polynomials\\

\hline
I. & $e^{-t^2};2t$ & $e^{\frac{t^2}{4}};\frac{t}{2}$ & $e^{2xt-t^2}=\sum_{n=0}^\infty H_n(x)\frac{t^n}{n!}$ & Hermite polynomials \\
    &                &                 &                                                 &   $H_n(x)$\cite{Andrew} \\

\hline
II. & $e^{-t^{m}};\nu t$ & $e^{(\frac {t}{\nu})^{m}};\frac{t}{\nu}$ & $\exp(\nu xt-t^{m})=\sum _{n=0}^{\infty }H_{n,m,\nu}(x)\frac {t^n}{n!}$ & Generalized Hermite  polynomials \\

& & &  &  $H_{n,m,\nu} (x)$ \cite{Lahiri} \\

\hline
III. & $(1-t)^{-1};$ & $(1-t)^{-1}$ & $\frac{1}{(1-t)}\exp{\Bigg(\frac{xt}{t-1}\Bigg)}=\sum _{n=0}^{\infty}L_{n}(x)t^n$ &  Laguerre  polynomials \\
& $\frac{t}{t-1}$ & $\frac{t}{t-1}$ &  & $n!L_{n}(x)$ \cite{Andrew}\\
\hline

IV. & $\frac{t}{1-t};~\ln \left(\frac{1+t}{1-t} \right)$ & $\frac{2}{e^{t}-1};~\frac{e^{t}-1}{e^{t}+1}$ & $\frac{t}{1-t}\left(\frac{1+t}{1-t}\right)^{x}=\sum _{n=0}^{\infty }\mathcal{P}_{n}(x)\frac{t^{n}}{n!}$ & $\mathcal{P}$idduck polynomials \\

& & & & $P_{n}(x)$ \cite{Boas,Erd3}\\

\hline
V. & $e^{\beta t};1-e^{t}$ & $(1-t)^{-\beta};$ & $\exp\Big(\beta t + x(1-e^{t})\Big)=\sum _{n=0}^{\infty }a_{n}^{(\beta)}(x)\frac{t^{n}}{n!}$ & Actuarial polynomials\\
&  & $\ln(1-t)$ & & $a_{n}^{(\beta)}(x)$ \cite{Boas}\\

\hline
VI. & $e^{-t}$; & $\exp\Big(a(e^{t}-1)\Big);$ & $e^{-t}\Big(1+\frac{t}{a}\Big)^{x}=\sum _{n=0}^{\infty }c{}_{n} (x\, ;a)\frac{t^{n}}{n!}$ & Poisson-Charlier  polynomials \\
& $\ln \Big(1+\frac{t}{a}\Big)$ & $a(e^{t}-1)$ & & $c{}_{n} (x\, ;a)$  \cite{Erd2,Jordan,Szego}\\

\hline
VII. & $\Big(1+(1+t)^{\lambda}\Big)^{-\mu};$ & $\Big(1+e^{\lambda t}\Big)^{\mu};$ & $\Big(1+(1+t)^{\lambda}\Big)^{-\mu}(1+t)^{x}$ & Peters polynomials \\
& $\ln (1+t)$ & $e^{t}-1$ & $=\sum _{n=0}^{\infty }s_{n} (x;\lambda,\mu)\frac{t^{n}}{n!}$ & $s_{n} (x;\lambda,\mu)$ \cite{Boas}\\

\hline
VIII.& $\frac{t}{\ln(1+t)};\ln (1+t)$ & $\frac{t}{e^{t}-1};e^{t}-1$ & $\frac{t}{\ln(1+t)}(1+t)^{x}=\sum _{n=0}^{\infty }b_{n} (x)\frac{t^{n}}{n!}$ &  Bernoulli polynomials \\
& & & &  of the second kind $b_{n} (x)$ \cite{Jordan}\\
\hline

IX.& $\frac{2}{2+t};\ln (1+t)$ & $\frac{1}{2}(1+e^{t});e^{t}-1$ & $\frac{2}{2+t}(1+t)^{x}=\sum _{n=0}^{\infty }r_{n} (x)\frac{t^{n}}{n!}$ &  Related polynomials\\
 & & & & $r_{n} (x)$ \cite{Jordan}\\

\hline
X.& $\frac{1}{\sqrt{1+t^{2}}};\arctan(t)$ & $sec t;\tan t$ & $\frac{1}{\sqrt{1+t^{2}}}\exp(x\arctan(t))=\sum _{n=0}^{\infty }R_{n}(x) \frac{t^{n} }{n!}$  &  Hahn polynomials\\
& & & & $R_{n} (x)$ \cite{Bender}\\
\hline

&&&&\\
XI.& $\left(1-4t\right)^{-{\frac{1}{2}}}$ & $\frac{1+t}{(1-t)^{a}};$ & $\left(1-4t\right)^{-{\tfrac{1}{2}} } \left(\frac{2}{1+\sqrt{1-4t} } \right)^{a-1}$ & Shively's pseudo-Laguerre \\
 & $\times \left(\frac{2}{1+\sqrt{1-4t} } \right)^{a-1};$ & $\frac{1}{4}-\frac{1}{4}\Bigg(\frac{1+t}{1-t}\Bigg)^{2}$ & $\times \exp\Bigg(\frac{-4xt}{\left(1+\sqrt{1-4t}~ \right)^{\, 2} } \Bigg)=\sum _{n=0}^{\infty }R_{n} (a,x) t^{n}$ & polynomials  $R_{n} (a,x)$ \cite{Rainville} \\
 & $\frac{-4t}{\left(1+\sqrt{1-4t}~ \right)^{\, 2}}$ & & &\\

\hline
\end{tabular}}}\\

{\tiny{
\begin{tabular}{|l|l|l|l|l|}
\hline
&&&&\\
XII.  & $\frac{1}{(1-t)^{1+\alpha+\beta}}$; & $\left(\frac{2}{1+\sqrt{1+2t}}\right)^{1+\alpha+\beta}$; &  $(1-t)^{-1-\alpha-\beta}{}_2F_1\left[
                                                                                                                                          \begin{array}{ccc}
                                                                                                                                            \frac{1+\alpha+\beta}{2}, & \frac{2+\alpha+\beta} {2};& \frac{2xt}{(1-t)^2} \\
                                                                                                                                             & 1+\alpha &  \\
                                                                                                                                          \end{array}
                                                                                                                                        \right]$
 &  The polynomials of \\
&$~\frac{2t}{(1-t)^2}$&   $~\frac{t}{1+t+\sqrt{1+2t}}$   &  $=\sum\limits_{n=0}^\infty J_n(x)\frac{t^n}{c_n}$  & the Jacobi case \cite{Roman1}\\
\hline

XIII.&$\frac{(1-t^2)}{(1+t)}$;   &   $\frac{1}{\sqrt{1-t^2}}$;& $\frac{1-t^2}{1-2xt+t^2}=\sum\limits_{n=0}^\infty (-1)^n t_n(x)t^n$& The polynomials of \\
&  $\frac{-2t}{(1+t)^2}$&  $\frac{-t}{1+\sqrt{1-t^2}}$  && the Chebyshev case \cite{Roman1}\\
\hline
\end{tabular}}}\\
\vspace{.35cm}

Now, we proceed to plot the graphs related to the 2ILP of order, $\alpha$ $L_n^{(\alpha)[2]}(x)$, 2ILP $L_n^{[2]}(x)$, 2IFF $(x)_n^{[2]}$ and 2IEP
$\phi_n^{[2]}(x)$ for $n=4$. For this, we need the
first few values of $L_n^{(\alpha)}(x)$, $L_n(x)$, $(x)_n$ and $\phi_n(x)$. We give the list of first five $L_n^{(\alpha)}(x)$,
$L_n(x)$, $(x)_n$ and $\phi_n(x)$ in Table 2.\\

\noindent
{\bf Table 2.~First five $L_n^{(\alpha)}(x)$, $L_n(x)$, $(x)_n$ and $\phi_n(x)$}.\\
\\
{\tiny{
\begin{tabular}{|l|l|l|l|l|l|}
  \hline
  $n$ & 0 & 1 & 2 & 3 & 4\\
  \hline

  $L_n^{(\alpha)}(x)$ &  1  &  $1+\alpha-x$  & $\frac{1}{2}(1+\alpha)(2+\alpha)$  & $\frac{1}{6}(1+\alpha)(2+\alpha)(3+\alpha)
  $ &$\frac{1}{24}(1+\alpha)(2+\alpha)(3+\alpha)(7+\alpha-x)-\frac{1}{8}
  (2+\alpha)(3+\alpha)(7+\alpha-x)x$ \\
  &&&$-(2+\alpha)x+\frac{1}{2}x^2$&$-\frac{1}{2}(2+\alpha)(3+\alpha)x $&$\frac{1}{8}(3+\alpha)(7+\alpha-x)x^2-\frac{1}{24}(7+\alpha-x)x^3-\frac{1}{8}(1+\alpha)(2+\alpha)(3+\alpha)$ \\
  &&& &$+\frac{1}{2}(3+\alpha)x^2-\frac{1}{6}x^3$& $+\frac{1}{24}(2+\alpha)(3+\alpha)x-\frac{1}{8}(3+\alpha)x$\\

\hline
  $L_n(x)$ & 1 & $1-x$ & $\frac{1}{2}(x^2-4x+2)$ & $\frac{1}{6}(-x^3+9x^2-18x+6)$ &$\frac{1}{24}(x^4-16x^3+72x^2-96x+24)$  \\

  \hline

  $(x)_n$  & $0$  & $x$&   $x^2-x$&   $x^3-3x^2+2x$&   $x^4-6x^3+11x^2-6x$\\
  \hline

  $\phi_n(x)$&   $1$&  $x$&  $x^2+x$&   $x^3+3x^2+x$&   $x^4+6x^3+7x^2+x$\\
  \hline
\end{tabular}}}\\
\vspace{.35cm}

In view of equations (4.10), (4.28), (4.37) and Table 2, we find the first few values of $L_n^{(3)[2]}(x)$, $L_n^{[2]}(x)$, $(x)_n^{[2]}$ and $\phi_n^{[2]}(x)$. We present the values of first five $L_n^{(3)[2]}(x)$, $L_n^{[2]}(x)$, $(x)_n^{[2]}$ and $\phi_n^{[2]}(x)$ in Table 3.\\

\noindent
{\bf Table 3.~First five $L_n^{(\alpha)[2]}(x)$, $L_n^{[2]}(x)$, $(x)_n^{[2]}$ and $\phi_n^{[2]}(x)$}.\\
\\
{\tiny{
\begin{tabular}{|l|l|l|l|l|l|}
  \hline
  $n$ & 0 & 1 & 2 & 3 & 4  \\

  \hline
   $\L_n^{(3)[2]}(x)$ &  $1$  & $x$  &  $\frac{x^2}{4}+\frac{5x}{2}-5$  & $\frac{11x^3}{36}-x^2-\frac{5x}{2}-30$  &  $\frac{x^4}{576}+\frac{7x^3}{48}
   +\frac{35x^2}{16}-\frac{35x}{24}-\frac{175}{8}$\\
  \hline

  $L_n^{[2]}(x)$ & $1$ &  $x$ & $\frac{x^2}{4}+x-\frac{1}{2}$  &  $\frac{x^3}{36}+\frac{x^2}{2}+\frac{x}{2}-\frac{13}{6}$& $\frac{x^4}{576}
  +\frac{5x^2}{8}-\frac{x}{6}-\frac{15}{24}$ \\
  \hline

  $(x)_n^{[2]}$ &  $0$& $x$  &  $x^2-2x$&  $x^3-6x^2+7x$  &   $x^4-12x^3+40x^2-35x$\\
  \hline

  $\phi_n^{[2]}(x)$  &  $1$ &  $x$  &  $x^2+2x$  &  $x^3+6x^2+5x$  &  $x^4+12x^3+32x^2+15x$\\
  \hline

\end{tabular}}}\\
\vspace{.35cm}

From Table 2, we get the following graphs:\\

\begin{figure}[htb]
\begin{center}

\epsfig{file=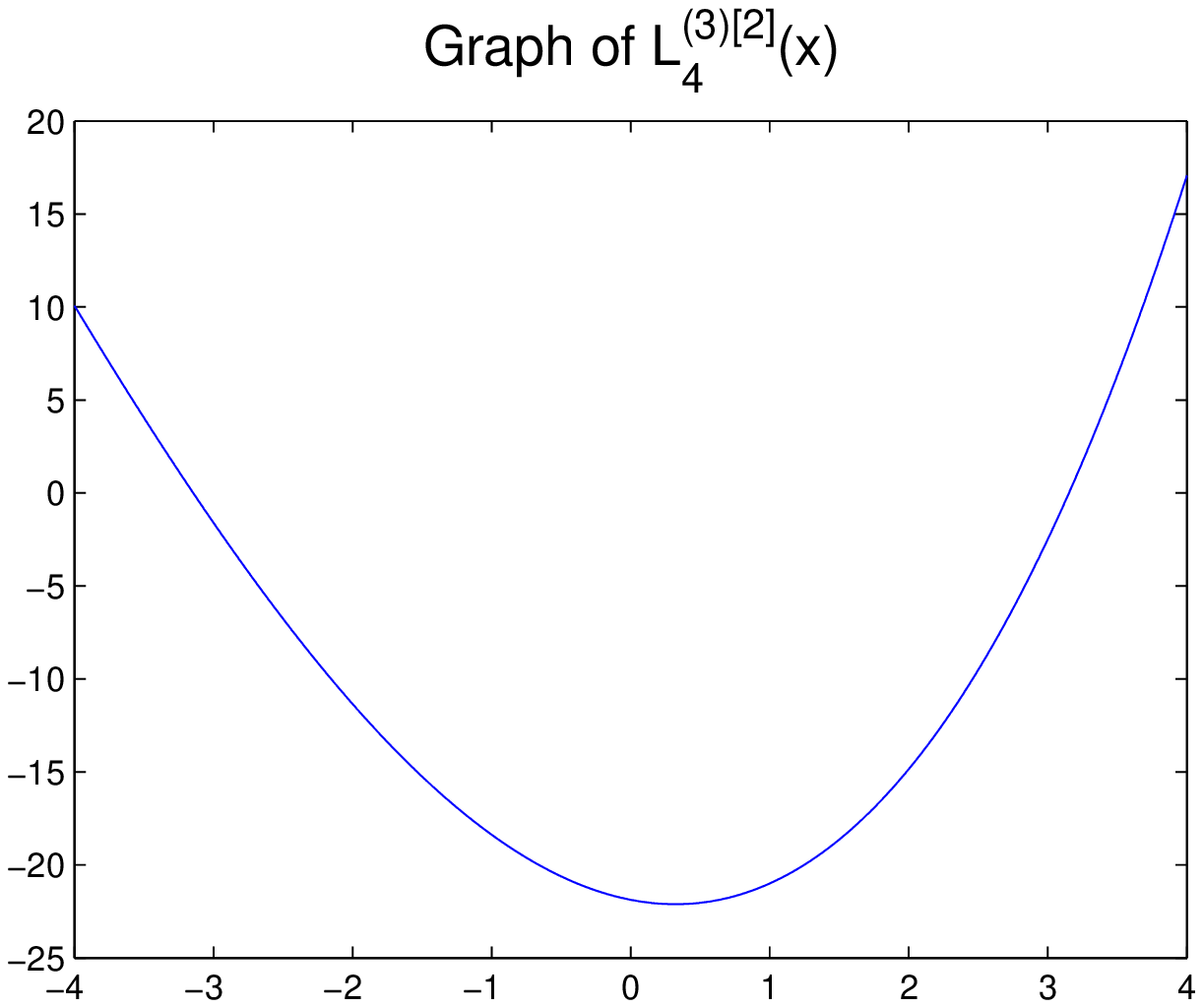, width=7cm}
\epsfig{file=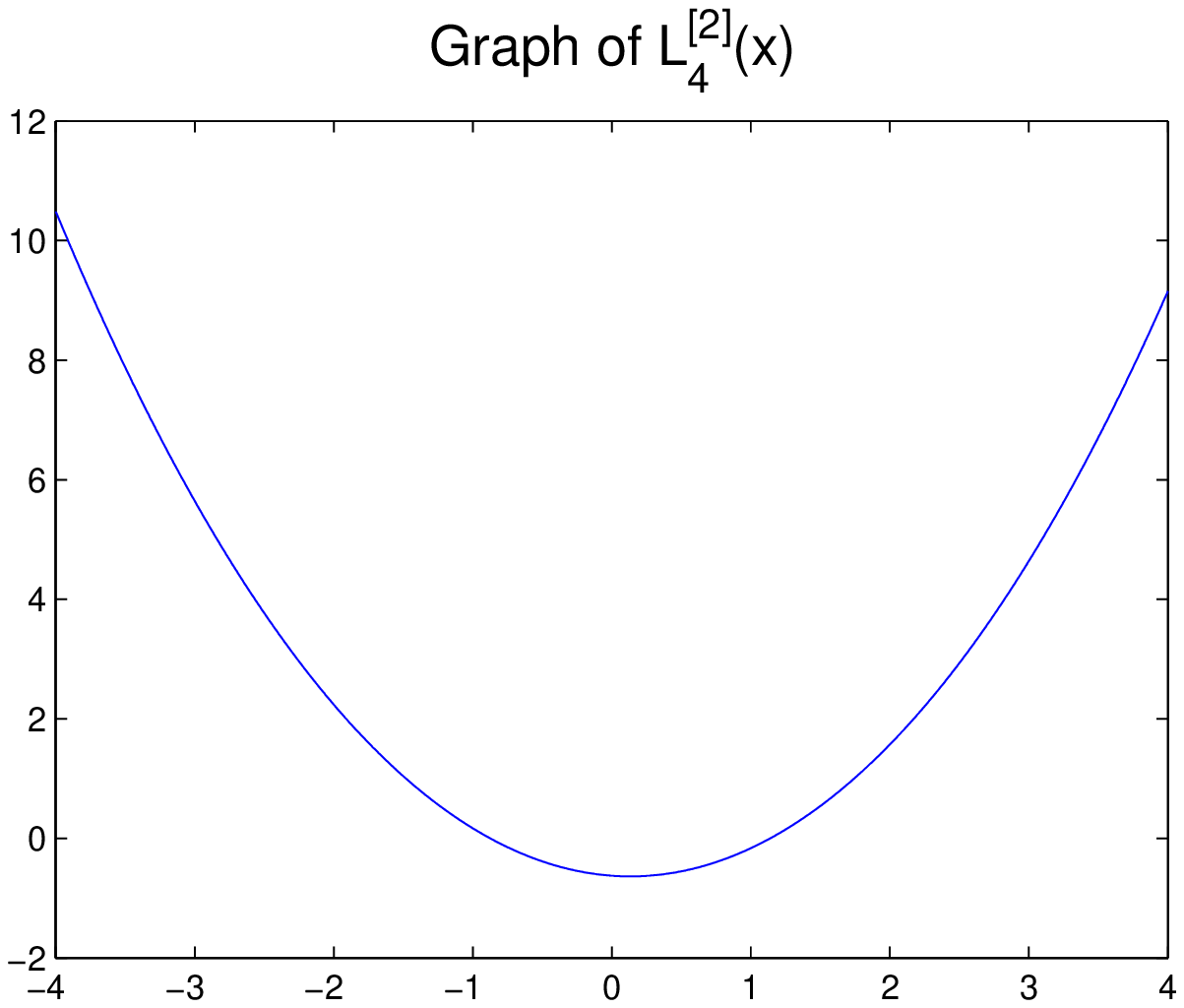, width=7cm}

\end{center}
\end{figure}

\begin{figure}[htb]
\begin{center}

\epsfig{file=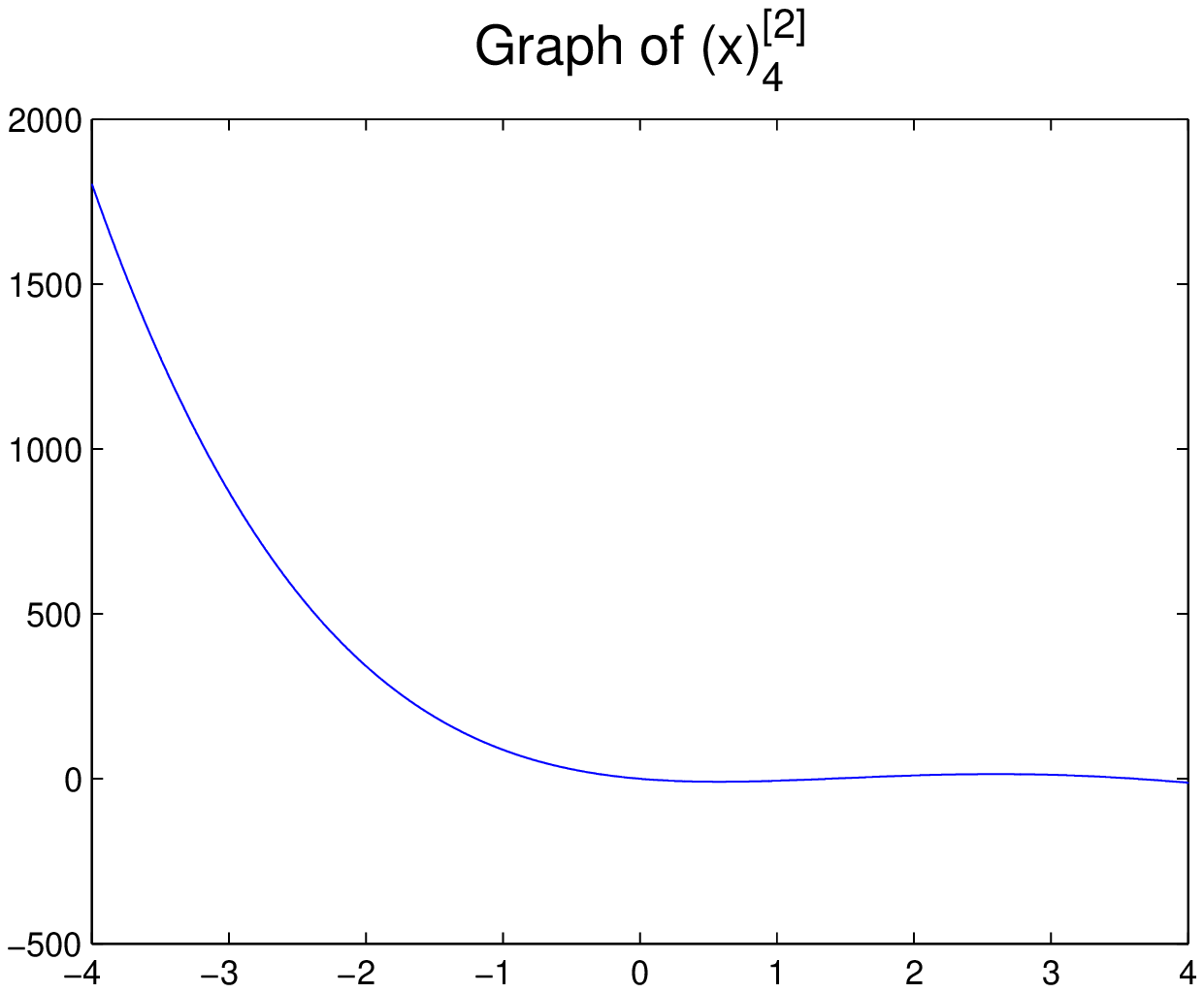, width=7cm}
\epsfig{file=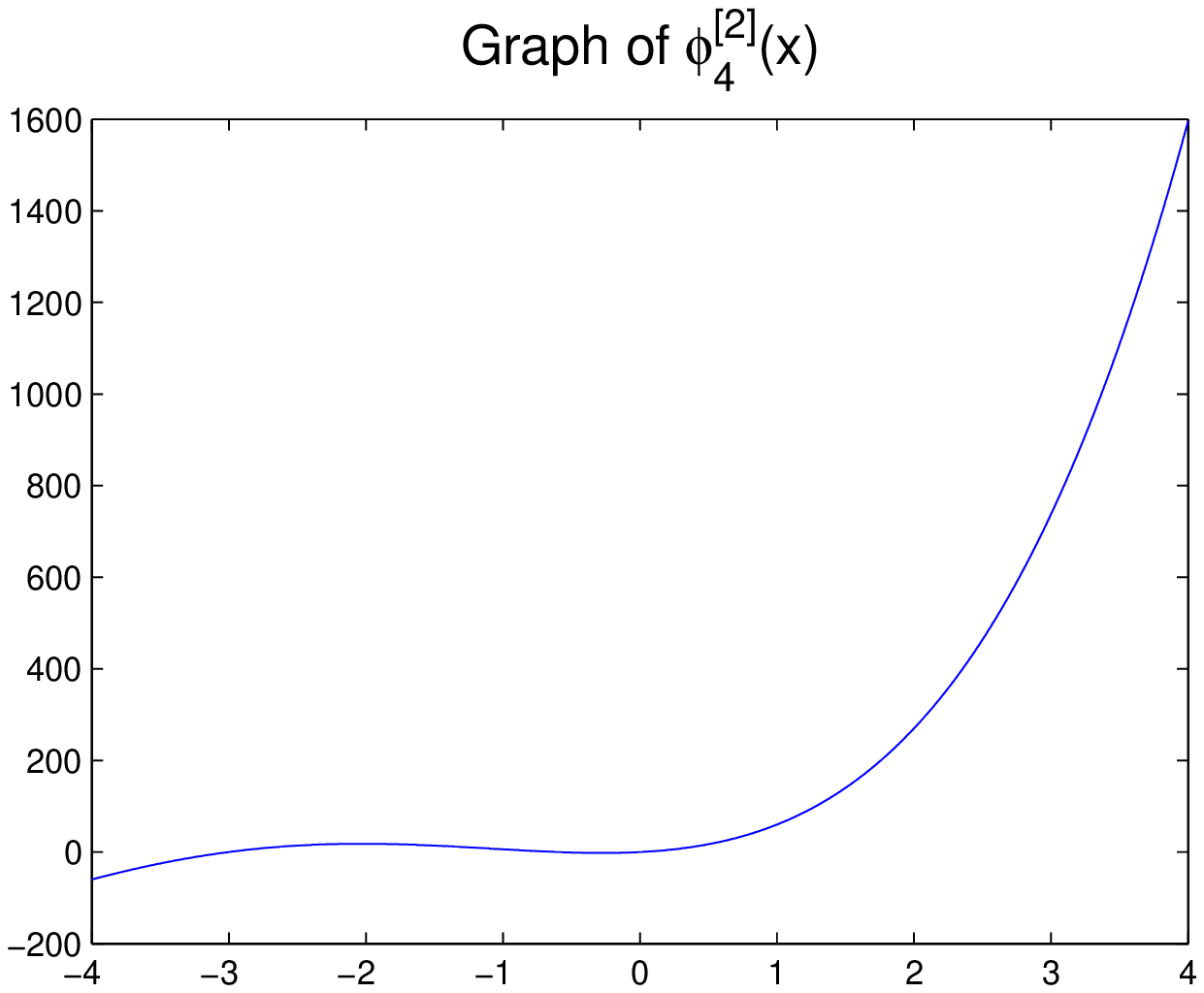, width=7cm}

\end{center}
\end{figure}

\end{document}